f

 \documentclass[preprint,review,12pt]{elsarticle}

 \usepackage{multirow}
 \usepackage{tabularx}
 \usepackage{rotating}

 \usepackage{graphicx}
 \usepackage{epsfig,subfigure,psfrag}

\usepackage{textcomp}
\usepackage{latexsym,amssymb,amsmath}

 \usepackage{lineno}

\begin{document}

\begin{frontmatter}



\title{A better method to determine the stability region of an
L-stable implicit-explicit Runge-Kutta scheme}


\author[a,b]{Shu-Chao Duan\corref{cor1}}\ead{s.duan@163.com}
\author{$et$ $al$}

\cortext[cor1]{Corresponding author}
\address[a]{Institute of Fluid Physics, CAEP, P. O. Box 919-105, Mianyang 621999, China}
\address[b]{Department of Modern Physics, University of Science and Technology of China, Hefei 230026, China}

\begin{abstract}
We propose a better method to determine the stability region of an L-stable
implicit-explicit Runge-Kutta scheme. This method always provides the correct result,
while other methods sometimes give wrong result. It is useful in
the analysis for implicit-explicit Runge-Kutta schemes.
\end{abstract}


\end{frontmatter}



\section{Introduction}\label{sec.introduction}

Many implicit-explicit (IMEX) Runge-Kutta (RK) schemes are widely used in the temporal integration of 
stiff-nonstiff mixed differential equations. The stability region is usually wanted in the analysis for an IMEX RK scheme. However, sometimes it is not easy to get it.

An IMEX RK scheme consists of applying an implicit RK scheme to the
stiff terms and an explicit scheme to the nonstiff terms. When applied to
$\partial_t U = R(U)+F(U)$, an $s$-stage IMEX RK scheme takes the following form \cite{Kennedy:2003, Pareschi:2005}:
\begin{eqnarray}
U^{(i)} &=& U^n + \Delta t \sum_{j=1}^s a_{ij} R(U^{(j)})
          + \Delta t \sum_{j=1}^{i-1} b_{ij}F(U^{(j)}),\label{eq.Ui}\\
U^{(n+1)} &=& U^n + \Delta t \sum_{i=1}^s w_i R(U^{(i)})
          + \Delta t \sum_{i=1}^{s} \omega_iF(U^{(i)}).\label{eq.Un1}
\end{eqnarray}
The matrices $A=(a_{ij})$ and $B=(b_{ij})$, where $b_{ij}=0$ for $j\ge i$, are coefficient matrices and correspond to the implicit and explicit part of the entire mixed scheme, respectively. The vectors $w = (w_1, ..., w_s)^T$ and $\omega=(\omega_1,...,\omega_s)^T$ 
are the weight vectors. So we can denote an IMEX scheme as $(A, w_1, B, \omega_2)$. Note that we placed the implicit and explicit parts in a different sequence compared with some papers because we have written them in the order indicated by the term IMEX itself: IM first, then EX.

The region of absolute stability $S_A$ of an IMEX scheme $(A,w_1,B, \omega_2)$ is defined as \cite{Kennedy:2003,Pareschi:2001}
\begin{equation}
S_A = \{(z_1, z_2)\in C^2 :|R(z_1, z_2)|\le 1\},
\end{equation}
where $z_1=\lambda_1 \Delta t, z_2 =\lambda_2 \Delta t$, $\lambda_1$ and $\lambda_2$ are the eigenvalues of the implicit and explicit part of the differential equation and $R(z_1,z_2)$ is the function of absolute stability of an IMEX scheme, defined as
\begin{equation}
R(z_1, z_2) = p(z_1, z_2)/q(z_1,z_2),
\end{equation}
where $p(z_1,z_2)=\texttt{Det}[I-z_1*A-z_2*B+z_1*e\otimes w_1 + z_2*e\otimes 
\omega_2]$, $q(z_1,z_2)=\texttt{Det}[I-z_1*A]$, where $I$ is the identity matrix and $e$ is the vector $\{1,1,...,1\}$.

$S_A$ is a region in $C^2 (\simeq R^4)$. It is not easy to visualize and comprehend the geometric structure of such a region. For L-stable schemes, instead of the region of absolute stability $S_A$, we use the L-stable region defined as \cite{Pareschi:2001}
\begin{equation}
S = \{ z_2 \in C:\texttt{Sup}_{z_1\in C^-}|R(z_1,z_2)|\le 1\},
\end{equation}
where $C^-$ is the complex left half-plane and $S$ is a region in $C (\simeq R^2)$.

The most straightforward method is to use the maximization function to realize the above definition. We call this the definition method for the later discussion. This method is too slow in our experience and sometimes provides an asymmetric result because it is difficult for the (global) maximization function to work well in the case of a high-order polynomial. Pareschi and Russo \cite{Pareschi:2001} proposed a second method that we call the parametric integration method, which is very fast. However, when the boundary of the stability region $S$ is not smooth, it provides the wrong result. Thus, we endeavored to provide a third method, which we call the algebra root method. It is much faster than the first method but slower than the second, and always provides the correct result. 

\section{New method}

To obtain this method, first define \cite{Pareschi:2001, Hairer:1987}
\begin{equation}
\begin{split}
f(y,\rho,\theta) &= |q(iy,-1+\rho e^{i\theta})|^2 - |p(iy, -1+\rho e^{i\theta})|^2, \\
g(y,\rho,\theta) &= \partial f / \partial y.
\end{split}
\end{equation}
Then, solve the following equations with a fixed $\theta$ in real domain:
\begin{equation}\label{eq.yRho}
\begin{split}
f(y,\rho) &= 0, \\
g(y,\rho) &= 0.
\end{split}
\end{equation}
This provides many roots. The smallest root of $\rho$ coupled with the prefixed $\theta$ consists of a point in the $(\rho,\theta)$ plane. A series of points produced by this procedure forms the boundary of the stability region $S$. This statement can be proved easily using Lemma 3.1 in \cite{Pareschi:2001} and other arguments, so it is omitted here for conciseness.

\section{Example}

\begin{figure}[htbp]
\centering
  \subfigure[]{\includegraphics[width=0.45\textwidth]{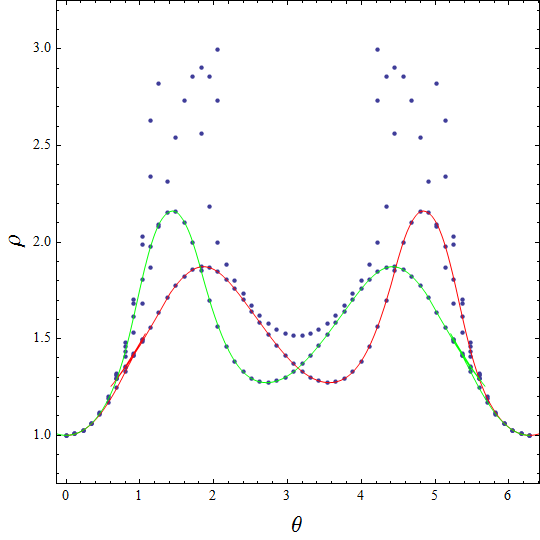}\label{fig.rhos}}
  \subfigure[]{\includegraphics[width=0.45\textwidth]{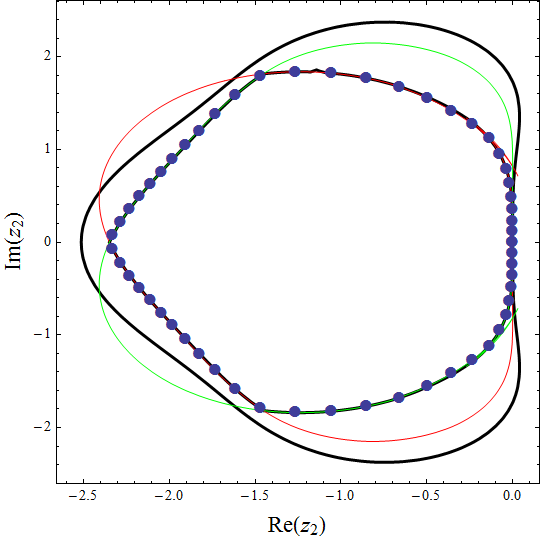}\label{fig.boundaries}}
\caption{Comparison of the three methods: black indicates the first method, red and green the second method, and dots the third method.}\label{fig.ComparionThree}
\end{figure}

We now provide an example to demonstrate these three methods. The example IMEX scheme is the ``additive RK.3.L.1'' referenced in the study by Liu and Zou \cite{Liu:2006}. In Fig. \ref{fig.rhos}, the dots are the roots of $\rho$ from the solutions of Eq. (\ref{eq.yRho}), where the abscissa is $\theta$. The red and green curves are integrated twice: once in one direction and again in the reverse direction, using the parametric integration method \cite{Pareschi:2001}. In Fig. \ref{fig.boundaries}, the two black curves are the boundaries of the stability region of the explicit part and the entire ``additive RK.3.L.1'' IMEX scheme, determined by the definition method. There is a small flaw in the inner curve, the boundaries of the stability region of the entire IMEX scheme, because of the inefficiency of the maximization function. The red and green curves are integrated by the parametric integration method as in Fig. \ref{fig.rhos}. The dots are obtained by the algebra root method. They coincide with the boundary determined by the definition method. We note that if there are more inflexions on the boundary, and hence the boundary is more complex, the parametric integration method does not work. Sometimes the integration diverges and does not follow the correct trace. We note additionally that the boundaries of the stability region of the IMEX scheme in the study by Liu and Zou \cite{Liu:2006} are slightly different from those in this paper. 

After the boundary is obtained, the area of the stability region can also be computed by first interpolating the boundary to a function $\rho(\theta)$ and then integrating along this trace: $\int_0^{2\pi}\frac{1}{2}\rho^2(\theta)d\theta$.

\section{Conclusion}

We propose a new method to determine the L-stable region of an L-stable IMEX RK scheme. 
This algebra root method is much faster than the
definition method but slower than the parametric integration method, and
always provides the correct result. It is useful in the analysis for IMEX RK
schemes.



\bibliographystyle{elsarticle-num}
\bibliography{<your-bib-database>}

\begin{thebibliography}{00}


\bibitem{Kennedy:2003}
C. A. Kennedy and M. H. Carpenter,
\newblock Appl. Numer. Math. 44 (2003) 139-181.

\bibitem{Pareschi:2005}
L. Pareschi and G. Russo,
\newblock J. Sci. Comput. 25 (2005) 129-155.

\bibitem{Pareschi:2001}
L. Pareschi and G. Russo,
\newblock Adv. Theory Comput. Math. 3 (2001) 269-289.

\bibitem{Hairer:1987}
E. Hairer and G. Wanner
\newblock Solving Ordinary Differential Equations, Vol.2 Stiff and Differential-algebraic Problems, Springer-Verlag, New York, 1987.

\bibitem{Liu:2006}
H. Liu and J. Zou, 
\newblock Journal of Computational and Applied Mathematics 190 (2006) 74-98.


 \end{thebibliography}



\end{document}